\documentclass[12pt, a4paper]{article}
\usepackage[T1]{fontenc}
\usepackage[utf8]{inputenc}
\usepackage{authblk}
\usepackage[pdftex]{graphics}
\usepackage{pdfpages}
\usepackage{amsmath,amssymb,amsthm}
\usepackage{pdfpages}
\usepackage[bottom]{footmisc}

\pdfoutput=1
\newcommand{\beq}{\begin{equation}}
\newcommand{\eeq}{\end{equation}}
\newcommand{\beqa}{\begin{eqnarray}}
\newcommand{\eeqa}{\end{eqnarray}}

\title{An extension of the Plancherel measure}
\author{Mikl\'os Arat\'o$^{1}${\thanks{\texttt{arato.miklos@renyi.mta.hu}}}, \;
Vill\H{o} Csisz\'ar$^{1}${\thanks{\texttt{csvillo@gmail.com}}}, \;
Bal\'azs Gerencs\'er$^{1,2}${\thanks{\texttt{gerencser.balazs@renyi.mta.hu}}}, \;
 Gy\"orgy Michaletzky$^{1}${\thanks{\texttt{michaletzky@caesar.elte.hu}}}, \;
L\'{\i}dia Rejt\H{o}$^{2,3}${\thanks{\texttt{rejto@udel.edu}}}, \;
 G\'abor Sz\'ekely$^{4}${\thanks{\texttt{gszekely@nsf.gov}}}, \;
 G\'abor Tusn\'ady$^{2}${\thanks{\texttt{tusnady.gabor@renyi.mta.hu}}} \;
 and  Katalin Varga$^{5}${\thanks{\texttt{vargaka@mnb.hu}}}
 }
\date{May 17, 2018}

\begin{document}

\maketitle

{\small{
\noindent $^1$ {E\"otv\"os Lor\'and University,
Department of Probability Theory and Statistics, Budapest, Hungary}\\
\noindent $^2$ {Alfr\'ed  R\'enyi
Mathematical Institute of the Hungarian Academy of Sciences,
Budapest,Hungary}\\
\noindent$^3$ {University  of Delaware, Department of Economics and Applied Statistics, Newark, DE, USA}\\
\noindent$^4$ {National Science Foundation, Washington D.C., USA}\\
\noindent$^5$ {Hungarian National Bank, Budapest, Hungary}
}}
\medskip

\bigskip \bigskip 
\quad \quad \quad \quad \quad \quad
\quad \quad \quad \quad \quad \quad
\quad \quad \quad \quad \quad \quad \quad
{\it Dedicated to the 80th birthday of Imre Csisz\'ar}

\medskip
{\small{
\noindent
{\bf Abstract}
Given a distribution in the unite square and having iid sample from it
the first question what a statistician might do to test the hypothesis that the sample is iid.
 For this purpose an extension of the Plancherel measure is introduced. 
Recent literature on asymptotic behavior of Plancherel measure is discussed with extension
 to the new set up. Models for random permutations are described and the power of different
tests is compared.
}}

\section{\bf {Introduction}}

\medskip
                                                            
Let $X(n), n=1,2,\ldots$ be iid uniform random variables on $(0,1)$. Step by step we define
the process 
$${\cal Y}(n) = \big(Y_{k,t}(n),t=1,\ldots,\lambda_k(n),k=1,\ldots,s(n)\big)$$
for all $n=1,2,\ldots$     
starting with $s(0)=0$. Here $s(n)$ is the number of levels, $\lambda_1(n),\ldots,\lambda_{s(n)}(n)$
are the Young numbers and $Y_{k,t}(n)$ are the new positions for the numbers
$X(1),\ldots,X(n).$
Note that the data do not change in course of the algorithm only we introduce a combinatorial structure
reflecting their value and order. The process ${\cal Y}={\cal Y}(n),n=1,2,\ldots$ is a dynamically structured ordered sample.

The random variables $s(1),s(2),\ldots$ generated by the $X(n)$ sequence will be monotone non decreasing
and for each level $k$ the positions $Y_{k,t}(n)$ have the  property of monotonicity;
$$
Y_{k,t}(n) < Y_{k,t+1}(n).
$$
For a demonstration, let us look  the first two steps of generating $Y$.
For $n=1$ we have $s(1)=1,\lambda_1(1)=1,$ and $Y_{1,1}(1)=X(1).$
\noindent
Next there are two possible cases:

--- if $X(1)<X(2)$ then $s(2)=1, \lambda_1(2)=2, Y_{1,1}(2)=X(1), Y_{1,2}(2)=X(2)$ and

--- $s(2)=2, \lambda_1(2)=\lambda_2(2)=1, Y_{1,1}(2)=X(2), Y_{2,1}(2)=X(1)$ otherwise.

The narrative of the algorithm is the following. For each $n$ the new number $X(n)$ sits down on the first level
to its place 
assigned according to its order of magnitude among the $X$-s. 
If there are any $Y$ on the first level
greater then  the newcomer then $X(n)$ pushes up the first $Y$ being larger 
then $X(n)$. On the next level the algorithm is repeated: the new number sits down on the place
ordered for it according to its magnitude and pushes up the first $Y$ on its level being larger.

For $i=1,\ldots,n$                                                                                                       
we denote by $\kappa_i(n)$ the level of the $i$th number in the ordered sample
after arriving the number $X(n),$ and by 
$\kappa(n)$ the                                                                                                                                                                                           
whole sequence. 

For example if $n=11$ and $X(1)=0.473, X(2)=0.117, X(3)=0.973, X(4)=0.832, X(5)=0.771, X(6)=0.032, X(7)=0.251,
X(8)=0.914, X(9) =0.343, X(10) =0.652, X(11)=0.574$ then $\kappa_1(11)=1, \kappa_2(11)=2,
\kappa_3(11)=1, \kappa_4(11)=1,
\kappa_5(11)=3, \kappa_6(11)=1, \kappa_7(11)=2, \kappa_8(11)=3, \kappa_9(11)=4, \kappa_{10}(11)=2, \kappa_{11}=5$.
Let us collect first the indices $i,$ where the value of $\kappa_i(11)$ equals one:


\begin{tabular}{cccc}                                                                                                                              
1 & 3 & 4 & 6 \\
\end{tabular}

\noindent
next where $\kappa_i(11)=2$:

\begin{tabular}{cccc}
2 & 7 & 10 \\
\end{tabular}

\noindent
next for $3$:

\begin{tabular}{cccc}
5 & 8 \\
\end{tabular}

\noindent
finally

\begin{tabular}{cccc}
9 \\
\end{tabular}

\noindent
and

\begin{tabular}{cccc}
11 \\
\end{tabular} 

\noindent                                                                                                        
for $4$ and $5$. The result is called standard Young tableau:  it is monotone increasing both in row and column.
These numbers are determined by the rank numbers of the $X(i)$-s:
$5,2,11,9,8,1,3,10,4,7,6$.
We call the $\kappa$ sequence generated by the rank numbers 
the {\bf level process} of a permutation.
Next table shows the evolution of the $\kappa$ sequences. Here, for the sake of clarity, we use
the inverse rank numbers in place of the sliding elements and we supplant with zeros the sleeping ones.

{\small{
\begin{tabular}{ccccccccccc}
1 & 2 & 3 & 4 & 5 & 6 & 7 & 8 & 9 & 10 & 11 \\
- & - & - & - & - & - & - & - & - & -  & -  \\
0 & 0 & 0 & 0 & 1 & 0 & 0 & 0 & 0 & 0  & 0  \\
- & - & - & - & - & - & - & - & - & -  & -  \\
0 & 0 & 0 & 0 & 2 & 0 & 0 & 0 & 0 & 0  & 0  \\
0 & 2 & 0 & 0 & 1 & 0 & 0 & 0 & 0 & 0  & 0  \\
- & - & - & - & - & - & - & - & - & -  & -  \\
0 & 0 & 0 & 0 & 2 & 0 & 0 & 0 & 0 & 0  & 0  \\
0 & 2 & 0 & 0 & 1 & 0 & 0 & 0 & 0 & 0  & 3  \\
- & - & - & - & - & - & - & - & - & -  & -  \\
0 & 0 & 0 & 0 & 2 & 0 & 0 & 0 & 0 & 0  & 4  \\
0 & 2 & 0 & 0 & 1 & 0 & 0 & 0 & 4 & 0  & 3  \\
- & - & - & - & - & - & - & - & - & -  & -  \\
0 & 0 & 0 & 0 & 0 & 0 & 0 & 0 & 0 & 0  & 5  \\
0 & 0 & 0 & 0 & 2 & 0 & 0 & 0 & 5 & 0  & 4  \\
0 & 2 & 0 & 0 & 1 & 0 & 0 & 5 & 4 & 0  & 3  \\
- & - & - & - & - & - & - & - & - & -  & -  \\
0 & 0 & 0 & 0 & 0 & 0 & 0 & 0 & 0 & 0  & 6  \\
0 & 0 & 0 & 0 & 6 & 0 & 0 & 0 & 0 & 0  & 5  \\
0 & 6 & 0 & 0 & 2 & 0 & 0 & 0 & 5 & 0  & 4  \\
6 & 2 & 0 & 0 & 1 & 0 & 0 & 5 & 4 & 0  & 3  \\
- & - & - & - & - & - & - & - & - & -  & -  \\
0 & 0 & 0 & 0 & 0 & 0 & 0 & 0 & 0 & 0  & 6  \\
0 & 0 & 0 & 0 & 6 & 0 & 0 & 0 & 7 & 0  & 5  \\
0 & 6 & 0 & 0 & 2 & 0 & 0 & 7 & 5 & 0  & 4  \\
6 & 2 & 7 & 0 & 1 & 0 & 0 & 5 & 4 & 0  & 3  \\
- & - & - & - & - & - & - & - & - & -  & -  \\
0 & 0 & 0 & 0 & 0 & 0 & 0 & 0 & 0 & 0  & 6  \\
0 & 0 & 0 & 0 & 6 & 0 & 0 & 0 & 7 & 0  & 5  \\
0 & 6 & 0 & 0 & 2 & 0 & 0 & 7 & 5 & 0  & 4  \\
6 & 2 & 7 & 0 & 1 & 0 & 0 & 5 & 4 & 8  & 3  \\
- & - & - & - & - & - & - & - & - & -  & -  \\
0 & 0 & 0 & 0 & 0 & 0 & 0 & 0 & 0 & 0  & 6  \\
0 & 0 & 0 & 0 & 6 & 0 & 0 & 0 & 7 & 0  & 5  \\
0 & 6 & 0 & 0 & 2 & 0 & 0 & 7 & 5 & 9  & 4  \\
6 & 2 & 7 & 9 & 1 & 0 & 0 & 5 & 4 & 8  & 3  \\
- & - & - & - & - & - & - & - & - & -  & -  \\
0 & 0 & 0 & 0 & 0 & 0 & 0 & 0 & 0 & 0  & 6  \\
0 & 0 & 0 & 0 & 6 & 0 & 0 & 0 & 7 & 0  & 5  \\
0 & 6 & 0 & 0 & 2 & 0 & 0 & 7 & 5 & 9  & 4  \\
6 & 2 & 7 & 9 & 1 & 0 &10 & 5 & 4 & 8  & 3  \\
- & - & - & - & - & - & - & - & - & -  & -  \\
0 & 0 & 0 & 0 & 0 & 0 & 0 & 0 & 0 & 0  &11  \\
0 & 0 & 0 & 0 & 0 & 0 & 0 & 0 &11 & 0  & 6  \\
0 & 0 & 0 & 0 & 6 & 0 & 0 &11 & 7 & 0  & 5  \\
0 & 6 & 0 & 0 & 2 & 0 &11 & 7 & 5 & 9  & 4  \\
6 & 2 & 7 & 9 & 1 &11 &10 & 5 & 4 & 8  & 3  
\end{tabular}
}}

The first element of the permutation is $5$ so we put the number $1$ in the fifth column.
Other elements remain zero.

The second element is $2$ (unfortunately) but this coincidence does not make any complication:
we put a $2$ in the second column and in the fifth, too.

The third is $11$ hence in the eleventh column appears the number $3$.

The fourth is $9$ and the number $4$ looks up in the ninth and the eleven column on the
appropriate level.

The fifth element is $8$: the number $5$ goes to the eighth column (and not $8$ in fifth).

And so on.
The number 
of rank numbers leading to this
$\kappa(11)$ sequence is
determined by the so called hook numbers of the tableau:

\medskip

\begin{tabular}{cccc}
8 & 5 & 3 & 1 \\
6 & 3 & 1 \\
4 & 1 \\
2 \\
1 \\
\end{tabular}

\medskip
                                                                                                                                                                              
 We know that $11!$ is divisible by the product of these numbers and the quotient  $2310$ is the number of
the rank numbers resulting to our $\kappa(11)$. 
Surprisingly  the number of
$\kappa(11)$s resulting the same $\lambda(11)$ sequence is $2310,$ too. 
The explanation is that the inverse of a permutation results in a $\kappa$ determining the same $\lambda$ as the $\kappa$
of the permutation and the two $\kappa$s determine uniquely the permutation.
This fact leads to the Plancherel distribution. The question is the joint density of
the numbers ${\cal Y}(n)$ and the limiting behavior of the process when $n$ goes to infinity.
The average divided by $\sqrt n$ and the standard deviation of $LP=-\ln(p)$ of the Plancherel distribution for some $n$ are the followings:

\medskip   

\begin{tabular}{ccc}
n  & ave & st.dev. \\
11 & 0.95 &0.9480  \\
$11^2$ & 1.46 & 1.8468 \\
$11^3$ & 1.72 & 3.7652 \\
$11^4$ & 1.82 & 6.7102 \\
$11^5$ & 1.85 & 12.4306 \\
$11^6$ & 1.86 & 19.6362 \\
\end{tabular}

\medskip

Standard deviation  is close to $0.57n^{0.25}$, cf \cite{Buf}, \cite{PM}, \cite{M}, and \cite{Ver}.
We use instead of $LP$ the statistic $H$ which is the sum of the logarithm of the hook numbers.
For $n=k*11^4, k=1,\ldots,11$ the process ${\cal Y}(n)$
is shown in Figure 1. The curves sample $Y_{k,t}$ only for some fixed $t$. 
Instead of $(0,1)$ the sample $X$ comes from uniform distribution in $(0,3000000)$ The figure shows the numbers $Y_{k,t}(n)$
according to their level: the first coordinate is $Y_{k,t}(n)$ and the second is $k$.

Figure 2 shows the  density estimation (red curve)
of $H$ for $n=11^5$, it is seemingly Gaussian (blue curve).

Figure 3 shows the $\lambda_k(n)$ numbers for $n=11^6$ offering the minimum for the product of
the hook numbers. The the minimum of $H$ is $11 859 210$ while its expected value is
$11 859 790$. The difference $580$ is relatively large concerning the standard deviation which
is $19.6362$.

The shape of the curves follow such a well
pronounced lines only when the corresponding Plancherel probability is close to its expected value.
It means that the sample elements are independent. 
We think that the process ${\cal Y}(n)$ offers an appropriate statistics for
testing independence with identical distribution.

In next table some parameters are given for $n=11^4$ on Ornstein-Uhlenbeck process.

\medskip

\begin{tabular}{cccc}
$\rho$    & ave   &    st.dev  \\
0.5     & 63008  &   3.00     \\
0.95    & 63017  &   9.23     \\
0.995   & 63129  & 101.35     \\
\end{tabular}

\medskip

Switching the distribution of $X(n)$ to standard exponential the habit of the process $\cal Y$ changes a bit
but the quantile transformation offers the appropriate chain between the two cases thus we use the same notation.
The distribution of $Y_{1,1}(n)$ is exponential with parameter $1/n$ hence we investigate the rescaled process
$${\cal Z}(n) = n{\cal Y}(n).$$ One can hope that the distribution of the process stabilizes as n goes to infinity; we call the limit
{\bf {extended Plancherel distribution.}}

\medskip

\section{\bf {The core process}}

\medskip

The statistics $S$ defined as $n!$ divided by the product of hook numbers counts all permutations belonging to the same level
process. The minimum of this number is one: there are two permutations being uniquely determined by their level process.
The maximum of the statistics is close to $\sqrt{n!}$. 
Let us sum up for a fixed $n$ for all integer $C$ the number of
permutations belonging to   $S$ 
 between $\exp(C-0.5)$ and $\exp(C+0.5)$.
 Then  for most $C$ the result is positive:
in such cases it has
to be larger than $\exp(2C).$ Surprisingly the distribution is concentrated in a narrow interval
containing the permutations what we call {\bf typical} (see \cite{BP}).

\medskip

\section{\bf {Conditional independence}}

\medskip

For fixed $n$ and $k<n$ the first $k$ element of a permutation $(p_i,i=1,2,\ldots,n)$ of the set $H=(1,2,\ldots,n)$ cuts the set $H$ in
two parts: the {\bf {initializing}} set $G_k=\{p_i,i=1,2,\ldots,k\}$ and its complement. 
For $1<k<n-1$ we say that  {\bf a random permutation can be cut properly}
by $k$ if $(p_i,i=1,2,\ldots,k)$ and $(p_i,i=k+1,\ldots,n)$ are conditionally independent on the condition
that $p$ is initialized by a given subset $G_k$ of $H$ with $k$ elements. We say that a random permutation is {\bf {proper (decomposable)}} if
it can be cut properly for all $1<k<n-1$.

If $n=11$ then there are $56$ Young tableaux and among them $27$ has larger Plancherel probability than $0.0076$. Their total probability is $0.951,$
thus we propose the test that accepts the uniform distribution of a random permutation if the corresponding Young tableau is an element
of this set with $27$ elements given in the next table.

\medskip

\begin{tabular}{ccccccc}
7 & 2 & 1 & 1 \\
6 & 4 & 1 \\
6 & 3 & 2 \\
6 & 3 & 1 & 1 \\
6 & 2 & 2 & 1 \\
6 & 2 & 1 & 1 & 1 \\
5 & 4 & 2 \\
5 & 4 & 1 & 1 \\
5 & 3 & 3 \\
5 & 3 & 2 & 1 \\
5 & 3 & 1 & 1 & 1 \\
5 & 2 & 2 & 2 \\
5 & 2 & 2 & 1 & 1 \\
5 & 2 & 1 & 1 & 1 & 1 \\
4 & 4 & 2 & 1 \\
4 & 4 & 1 & 1 & 1 \\
4 & 3 & 3 & 1 \\
4 & 3 & 2 & 2 \\
4 & 3 & 2 & 1 & 1 \\
4 & 3 & 1 & 1 & 1 & 1 \\
4 & 2 & 2 & 2 & 1 \\
4 & 2 & 2 & 1 & 1 & 1 \\
4 & 2 & 1 & 1 & 1 & 1 & 1 \\
3 & 3 & 3 & 1 & 1 \\
3 & 3 & 2 & 2 & 1 \\
3 & 3 & 2 & 1 & 1 & 1 \\
3 & 2 & 2 & 2 & 1 & 1 \\
\end{tabular}

\medskip

The most probable tableau is $5,3,2,1,$ its Plancherel probability if $0.1337.$ If we choose a random distribution
having uniform distribution on the set of permutation having a corresponding Young tableau equal with this one
the divergence projection of this distribution on the set of distributions which may be properly cut the resulting
distribution has a divergence from the uniform one $7.5$.

Proper random permutations have a simple parametrization: for all subset $G$ of the set $H=\{1,\ldots,n\}$
there is a probability distribution. Defining the elements of the random permutation step by step we use the
distribution on the subset of the remaining elements.

We call the random permutation {\bf {double proper}}
 if both the random permutation and its inverse are
proper. If all permutations have positive probability the logarithms of the probabilities are elements of
a subset with dimension $\sum^{n-1}_{i=1} i^2$ (\cite{Cs}). The space itself has a rather sophisticated structure but
there is a subspace with simple characteristics:
$$
\log P(\pi) = \sum^n_{i=1} a_{i,\pi(i)},
$$
where $A=a_{i,k}$ is a real valued $n\times n$ matrix. We call the model {\bf checkerboard}.
In case we have {\bf only one permutation} even the inference
of such a simple model is out of question without further structural assumptions on the matrix $A$. The graph of
the two dimensional points $(i,\pi(i),i=1,\ldots,n)$ is a simple pictorial statistics of a permutation. If we
can figure out any structure in the graph the permutation is certainly not random and there is a hope to fit
a matrix $A$ with relatively small number of degree of freedom.

Suppose we have a two-dimensional random variable $(X,U)$, then we can construct the same processes $V_i,\tau_i$
from iid sequence $U_i$ as we have done from $X_i$ and the pair of processes may reflect the joint distribution
of $(X,U)$. There is an other possibility: let $\pi_i(n)$ be the permutation moving the ordered sample of
$X_1,\ldots,X_n$ to the ordered sample of $U_1,\ldots,U_n$ than the level process of this permutation
reflects directly the joint distribution of $(X,U)$. If $X$ and $U$ are independent, then the distribution of $\pi$ 
is uniform and the distribution of the level process is Plancherel. In the general case the distribution 
resembles to the checkerboard one. If the pair $(X,U)$ has two dimensional normal distribution then the
distribution of the level process depends only on the correlation of $X$ and $U$. We call the process
{\bf Gaussian level process} 
(see \cite{Co} and \cite{S}).

\medskip

\section{\bf {Complexity}}

\medskip

There are many statistics measuring different aspects of permutations. We can invert any permutation into some
iid sequence reordering its ordered sample according to the given permutation. Applying the method we can resample
a single permutation multiplying it with iid uniformly distributed random permutation and seeking any deviance
between the original permutation and its random descendants. One possible method is compression: if there are any
short description of the given permutation it is not highly
complex (\cite{A}, \cite{C}, \cite{HaNa} and \cite{KMSW}).

A possible way is the extension of Lov\'asz's graphons, see \cite{HKMS}, \cite{KKRW} and \cite{FMN}.
The idea is to order to all subset $A$ with $k$ elements of a permutation $\pi(i)$ of the first $n$ integer the
rank numbers corresponding to the numbers $\pi(i), i\in A.$

\medskip

\section{\bf {Dynamics}}

\medskip

One face of a permutation is that it represents a dynamics: $\pi_i$ means that {\it number $i$ moves to number $\pi_i$}
what makes sense even in case $\pi_i=i.$ Using cyclically the neighboring pairs $\pi_i,\pi_{i+1}$ we can define a
new permutation $\delta$ as
$$
\delta(\pi_i) = \pi_t, \quad {\text where} \quad t=i+1 \quad {\text if} \quad i<n \quad {\text and} \quad  t=1 \quad {\text otherwise}.
$$

\noindent
The attractors of this dynamics may reflect the complexity of $\pi$.

\medskip

\section{\bf {Testing IID property}}

\medskip

In nonparapetric statistics the sample elements are supplanted by their rank numbers
\cite{WW}, \cite{ChoW} and \cite{FP} and the typical resampling method is the use of
random permutations on some parts of the sample. We are trying to bridge this traditional
field of statistics with Plancherel measure.

\medskip

\section{\bf {Variations}}

\medskip

Victor Reiner, Franco Saliola and Volkmar Welker initiate the use of $k$ elements subsets of a permutation
{\it without} supplanting the elements with their rank numbers \cite{RSW}. They conjecture that  all the
square of the singular values of the matrix joining these parts with the original permutations are integers.

\medskip

\section{\bf {Power studies}}

\medskip

For arbitrary real $t$ let us define the distribution $P_t$ by
$$
P_t = c_t\exp(tLP),
$$
where $c_t$ is an appropriate scaling factor (of course $c_0=1$).
First question is the relation of our test with other ones concerning this exponential family.
It would be interesting to compare the power of hook-number-product test against this
exponential family with the power against the Gaussian level alternative (see \cite {LC}).

\medskip

\section{\bf {Where does the information come from?}}

\medskip

Once upon a time, fifty years ago, there was a conference in Debrecen. At this conference
Shinzo Watanabe gave a lecture with the title we are using here again
(\cite{W}).
The truth is that our world was originally full of information. Ir was the Demiurge who filtered
out superfluous parts of this information and imposed pre-existing Forms on the chaotic material. 
At the above mentioned conference Imre Csisz\'ar, Gyula Katona an G\'abor Tusn\'ady proved  the conservation
of entropy. Now we think that the controlled loss of information makes the trick: when we reduce
a permutation to its level process we condensate the information of the permutation in a proper way.

\medskip

\section{\bf {Laudation}}

\medskip

Many thanks to Imre Csisz\'ar for providing us the privilege to receive and benefit from many wise advice
given us by him throughout
 our  life and, in particular, for insightfully encouraging us to write on the present theme
a paper that we are pleased to dedicate to him on the occasion of his 80th birthday, and to wish him happy returns in an
arbitrarily long life in good health and spirits.

\includepdf[pages=-]{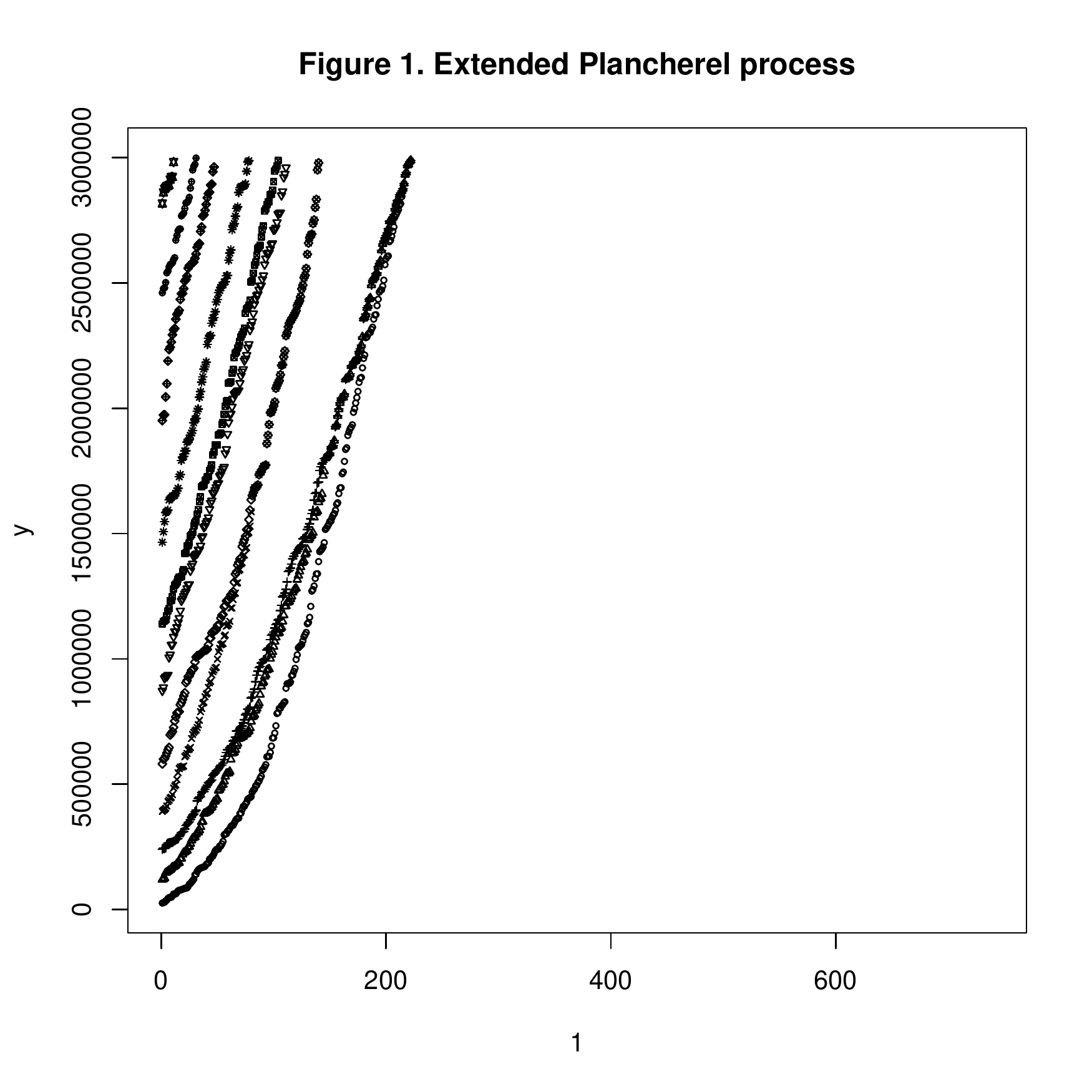}

\includepdf[pages=-]{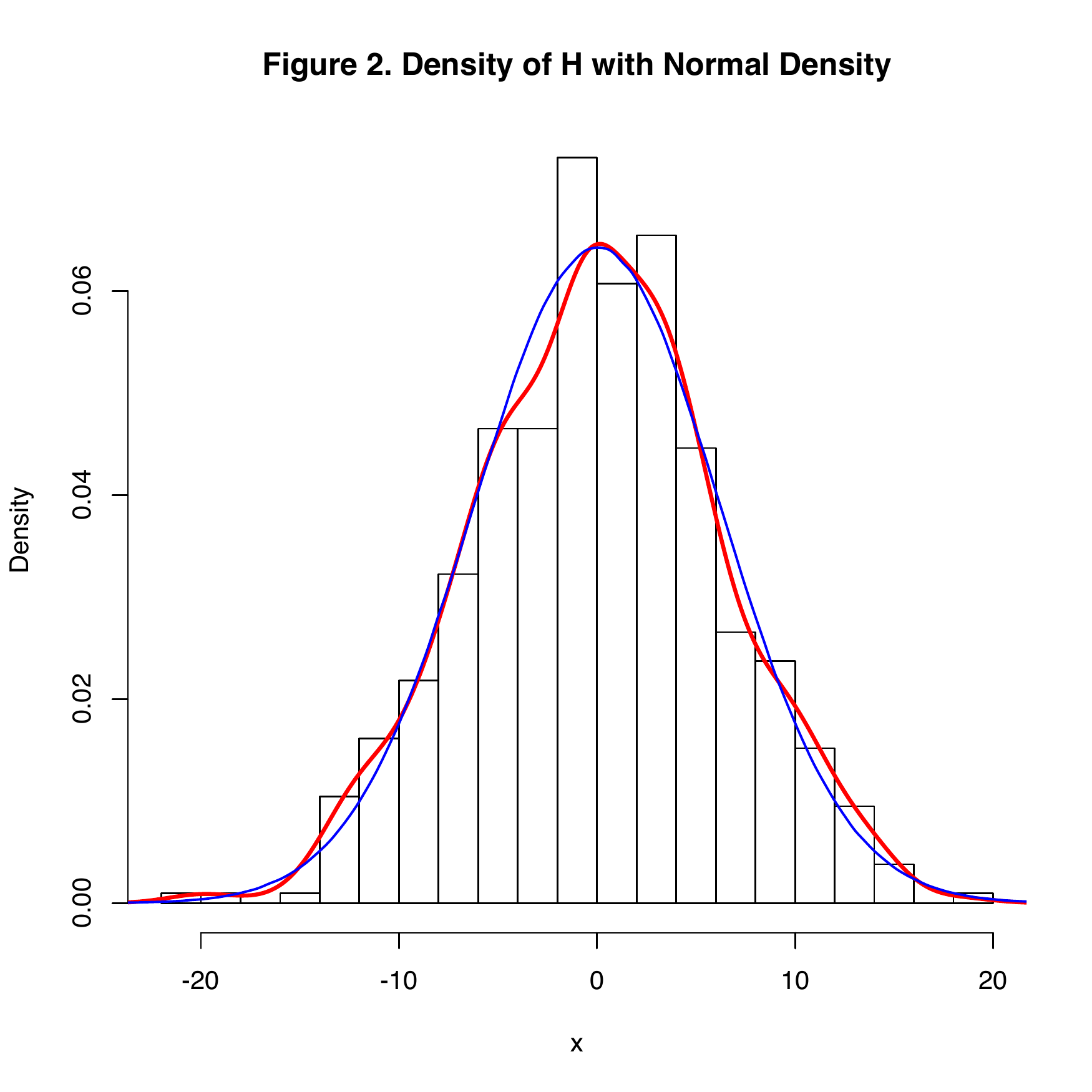}

\includepdf[pages=-]{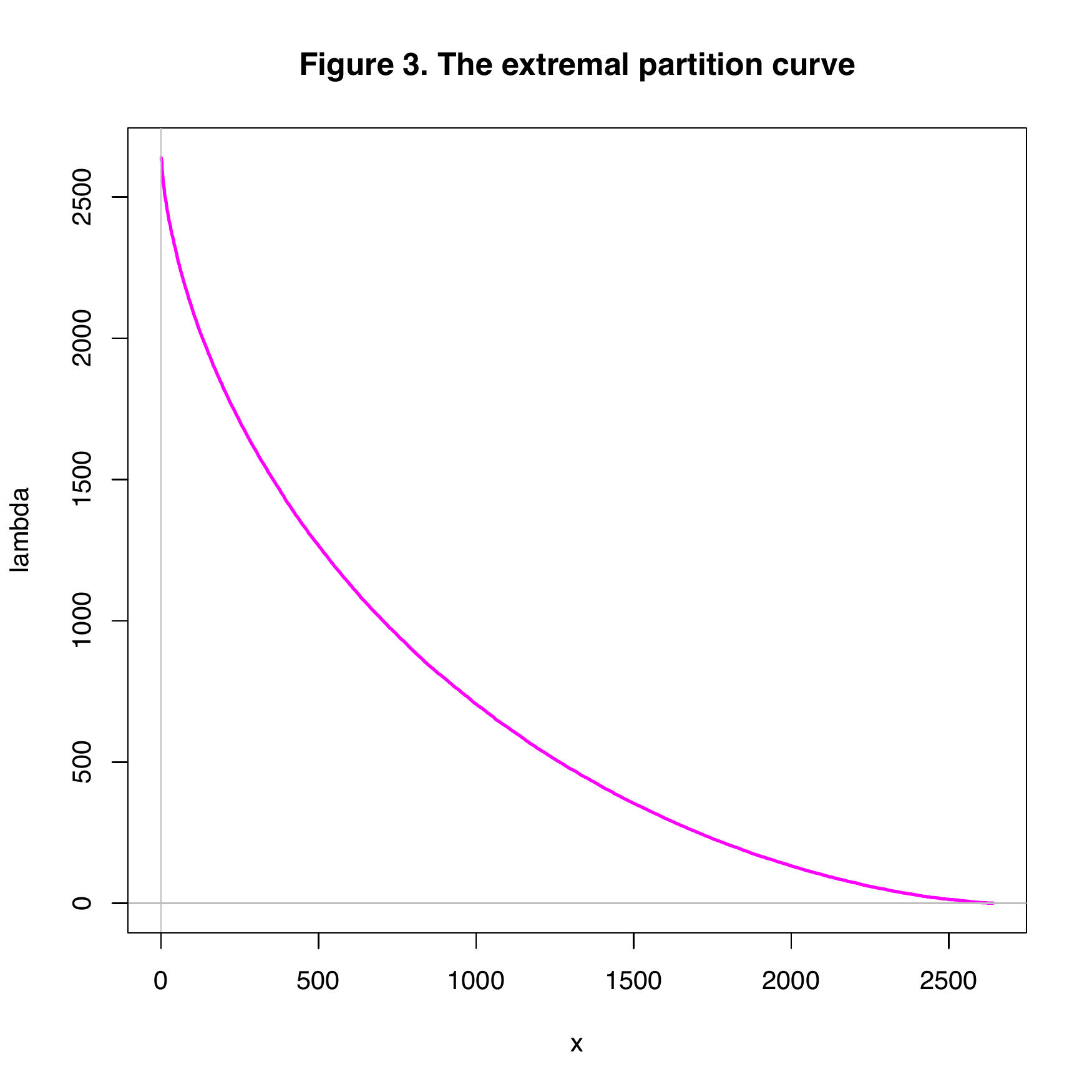}

\end{document}